\documentclass{elsart3}
\usepackage{amsmath}
\usepackage{graphicx}
\usepackage{amssymb}
\usepackage{natbib}
\newcommand{\Mc}{\mathcal{M}}
\newcommand{\Mb}{\mathbf{M}}
\newcommand{\Ab}{\mathbf{A}}
\newcommand{\x}{\mathbf{x}}
\newcommand{\y}{\mathbf{y}}
\newcommand{\z}{\mathbf{z}}

\newcommand{\Mn}{\mathcal{M}_n}
\newcommand{\Sc}{\mathcal{S}}

\newcommand{\Esym}{\text{E}}
\newcommand{\E}[1]{\Esym\left[#1\right]}
\newcommand{\Probsym}{\mathbb{P}}
\newcommand{\Prob}[1]{\Probsym\left[#1\right]}

\begin{document}
\begin{frontmatter}
\bibliographystyle{elsart-num}
\title{A set theoretic framework for enumerating matches in surveys
and its application to reducing inaccuracies in vehicle roadside surveys}
\author{Richard G. Clegg}
\ead{richard@richardclegg.org}
\address{Department of Mathematics, University of York,\\
York, YO10 5DD, United Kingdom}

\begin{abstract}
This paper describes a framework for analysing matches in multiple data
sets.  
The framework described is quite general and can be applied to a variety of 
problems where matches are to be found in data surveyed at
a number of locations (or at a single location over a number of days).  As
an example, the framework is applied to the problem of false matches
in licence plate survey data.  The specific problem addressed is that
of estimating how many vehicles were genuinely sighted at every one
of a number of survey points when there is a possibility of accidentally
confusing two vehicles due to the nature of the survey undertaken.

In this paper, a method for representing the possible {\em types of
match} is outlined using set theory.  The phrase {\em types of match}
will be defined and formalised in this paper.  A method for enumerating
$\Mn$, the set of all types of match over $n$ survey sites,
is described.  The method is applied to the problem of correcting
survey data for false matches using a simple probabalistic method.  An
algorithm is developed for correcting false matches over multiple
survey sites and its use is demonstrated with simulation results.
\end{abstract}


\end{frontmatter}

\section{Introduction}

In the analysis of roadside survey data, it is often desirable to
analyse matches between several data sets simultaneously.
For example, we might wish to answer questions
of the general type
``How many drivers are seen at point A, point B and point C?'' or
``How many vehicles are seen on all five survey days?''  This paper
attempts to create a general framework for the analysis of matching 
between data from more than two surveys. The framework is then
applied to the specific case of false matching in partial licence
plate surveys (that is non-matches which are mistaken for matches
because only part of the licence plate is observed).  It should be 
stressed throughout that the framework outlined is applicable to
any data series where matches are sought between two or more
distinct data sets.  While the work is placed in the context of
licence plate surveys (and further in the context of licence plate
surveys using a specific type of British licence plate) the results
are much more general than this.

Licence plate surveys are commonly used in the study of traffic systems,
particularly when measurements of the same vehicle are
required more than
one point (for example, calculating travel times or the routes of
vehicles).  Although automated techniques are becoming more common (GPS,
toll-tags and automated recognition cameras) the manual
licence plate survey remains an important tool for the road transport
engineer.  If a road with a high volume of traffic is being surveyed
then it is often the case that only part of the licence plate is recorded.
When this is the case, the possibility of spurious matches occurs.  To
take an example, standard British licence plates used to be of the
following
form: single letter, three digits, three letters: e.g. {\tt A123BCD}.
This form will be used throughout the paper, however, it must be 
stressed that this method would work with partial observations
of any type given the assumptions stated later.  If a surveyor
only recorded the first letter and three digits, then a vehicle {\tt A123ABC}
would not be distinguished from a vehicle {\tt A123XYZ} since the
disambiguating information (the final three letters) would not be recorded.

While the chances of such a false match are low, quite often the combinatorics
of the problem means that the actual recorded number of false matches remains
high.  To mathematicians, this is familiar as the celebrated {\em Birthday
Paradox}.  The Birthday Paradox asks the question ``How many people must
we have in a room before we might expect that two 
share the same birthday?''  Intuitively,
we might expect this to be quite a high number (since it is unlikely that
any two people share a birthday).  However, the number of pairs of people
in a room goes up with the square of the number of people in the room
$(n^2-n)/2$.  If we made the assumption that the chance of two randomly
selected people sharing a birthday is one in 365 then we only need twenty three
people in the room before it becomes likely (probability above 50\%) that two
will share a birthday.  Combinations in multiple point surveys work similarly.
If we had two survey sites, each with one thousand observations then
this is one million pairs of observations.  If the chances of a false match
in a given pair are only one in a ten thousand, we will still get (on average)
one hundred false matches.  This could well be larger than the actual number
of genuine matches in the data set and will certainly be a significant bias.

This paper attempts to provide a sound theoretic backing (using the well-known
framework of set theory) to matching problems across multiple data sites.
In section two, a general background of matching problems in licence plate data
is given to put the problem into context within the transport field.
In section three, the concept of {\em types of match} is formalised using the
standard set theoretic concept of an equivalence class.
In section four, a simple method is given for constructing the set $\Mn$, 
the set of every possible type of match across $n$ survey sites.  In section
five, partial ordering is introduced to apply the problem to false matches
due to incomplete observations.
In section six, an algorithm is given for correcting false matches using 
the framework developed in sections three, four and five.
Finally, in section seven, computational results are
given on artificially generated survey data.  The work in this paper
can be found in a much expanded form in \cite[Chapter Four]{clegg2004} and
an example of the method being used on real road traffic 
data is found in \cite[Chapter Five]{clegg2004}.  The set theory used
in this paper is extremely simple (just the concepts of equivalence class
and partial order are necessary) and would be covered in any standard
text on the subject, for example \cite{halmos70}.

\section{The false match problem in licence plate data at multiple sites}

Throughout this paper, the 
examples are given using an old form of British licence plate --- it should be
stressed that this is not necessary for this framework and is done purely for
the sake of example.  The work described here assumes nothing about the nature
of the individuals being observed other than the restrictions described in
Definition \ref{defn:pi}.  Similarly, when the phrase observation sites is used 
throughout this paper, this can mean either geographically distinct observation sites or  
a single geographical location observed for a number of days or any
combination of times and locations.  (In the work which motivated this research, the 
experimenters were interested in finding vehicles which travelled between three distinct
geographical locations on two consecutive days.  This would count as six observation sites
in the terminology used here.)  Note that no time information is used here although time
information is often available for such surveys.  It is hoped that a future improvement
to this method will make use of time information about observations to reduce uncertainties.

It is often the case that on-street traffic surveys collect partial vehicle
licence plate information. [The reason for
collecting partial rather than full licence plate information is that the
recording and transcription of the data is often done manually and time
constraints would preclude recording a full plate.]  This information can
then be used to reconstruct travel times and to infer route information about
drivers.  In partial plate data, however, problems can occur from
{\em false matches} as discussed above.  Of course, false matches
could also occur through recording or transcription errors.  While this
paper will not discuss these problems, it is in principle possible to 
extend this framework to cover recording and transcription errors.

In the case of two survey sites and no recording or transcription errors
the situation is relatively clear.
If our data shows that a match occurs between two observations (one from each
site) then, this must mean that either the same vehicle has been observed
at both, or that two different vehicles have been observed which happened to
have the same partial licence plate.
At multiple sites the situation is much more complex.  An apparent match
at four survey points may be any of the following: 
a true match (the same vehicle seen at all four
points); a different vehicle at each of the four points which (by coincidence)
have the same partial plate; a vehicle at survey point one and two which
has the same partial plate as a second vehicle at survey points three and
four; or any other of fifteen total possibilities.  The problem
becomes more difficult as the number of sites increases.  Indeed it is
not immediately clear how to enumerate the number of ways in which a match
as described above can occur over multiple data sites.  This issue is not a
trivial one.  In real licence surveys, the number of false matches is often
greater than the number of true matches.  In  
\cite[Chapter 5]{clegg2004} two survey sites with a flow of approximately 
one thousand vehicles at each were found to have ninety observed matches
between vehicles despite the fact that (given the positioning of the sites) it
would be extremely unlikely for any drivers at all to travel between them.

A number of researchers have approached the false matching problem for licence
plates. 
An early approach for two sites is given by \cite{hau1} which uses
a simple probabalistic correction. Several methods are described in
\cite{mah1}
including the possibility of two point matches between vehicles 
observed at
pairs of sites selected from several survey sites (for example entering and
leaving a cross-roads).  A graphical procedure
for visualising matches based upon journey time
between two sites is given by \cite{watling1988}.
Methods in this paper are useful for any analysis of data in which 
time between observations is a factor.  Further refinements for
site pairs, including a maximum likelihood method based upon assumptions
about travel time distribution are given in
\cite{watling1992} and \cite{watling1994}.
However, all of these methods concentrate on matches between pairs of sites
and the majority of them also assume that journey time information can be
used to aid in finding false matches,
which is not the case if, for example, we are interested in correcting
false matches at the same site over different days.  The
method described in this paper concentrates on matches between observations at
more than two sites, particularly where journey time information is not available or
cannot be used.

It should be emphasised again that, while this work is presented within the context of
licence plate surveys (indeed within the context of licence plate surveys on a
specific type of British licence plate) the results presented are extremely general.
These results would be applicable to any type of survey data where individuals are sought
in more than two data sets and where a possibility of confusion between observations
of individuals exists.  Applications for this technique are being sought in other
areas such as DNA matching and suggestions for suitable data sets would be welcomed
by the author.

\section{Equivalence classes for representing types of match}

In this section, notation is given, with examples, to describe a mathematical 
framework for investigating matches in multiple data sets.  For the convenience
of the reader the notation used throughout this paper
is gathered here for reference and defined as it
occurs throughout the paper.  In general bold lower case $\x$ is used to indicate 
a tuple (ordered set).  Upper case $M$ is used to indicate a set and bold upper
case $\Mb$ is used to indicate a set of sets.  Caligraphic lettering $\Sc$ is used
to indicate higher order entities such as sets of tuples or sets of sets of sets.

The following specific notation is used.
\begin{itemize}
\item $n$ --- the number of sites under investigation.
\item $\# M$ --- the number of members of set $M$.
\item $S_i$ --- the set of observations at site $i$.  See Definition \ref{defn:Si}.
\item $\y$ --- a tuple of observations, one from each site.  See Definition \ref{defn:y}
\item $\Sc$ --- the set of all possible tuples of observations.  See Definition \ref{defn:Sc}
\item $\Mb = \{M_1, M_2, \dots, M_m\}$ --- a {\em type of match}.  See Definition \ref{defn:Mb}.
\item $\Mn$ --- the set of all types of match for $n$ sites.  See Definition \ref{defn:Mn}.
\item $C(\y)$ --- the {\em type of match} of a tuple of observations $\y$.  See Definition
\ref{defn:Cy}.
\item $\Ab_n$ --- the set of sets $\{ \{1,2, \dots, n \} \}$ representing the same
observation across all sites.  See Definition \ref{defn:taun}.
\item $\y^*$ --- the tuple of {\em partial observations} from the tuple $\y$.  See Definition
\ref{defn:partial}.
\item $\Sc^*$ --- the set of all such {\em partial observations}.  See Definition 
\ref{defn:partial}.
\item $x(\y,\Mb)$ --- the {\em exact matching function} for the tuple $\y$.  See Definition 
\ref{defn:xyM}.
\item $X(\Mb)$ --- the exact matching count for the set $\Sc$.  See Definition 
\ref{defn:XM}.
\item $r(\y,\Mb)$ ---  the {\em relaxed matching function} for the tuple $\y$.  See Definition 
\ref{defn:ryM}.
\item $R(\Mb)$ --- the relaxed matching count for the set $\Sc$.  See Definition 
\ref{defn:RM}.
\item $T(M)$ --- the number of observations which are the same across 
all sites in the set $M$.  See Definition \ref{defn:TM}.
\item $p(i)$ --- the probability that $i$ distinct individuals, different in a full
observation, are the same in a partial observation.  See Definition \ref{defn:pi}.
\end{itemize}

\begin{defn}
Let $n$ be the number of observation sites and let $S_i$ be the set of observations at the
$i$th such site.
\label{defn:Si}
\end{defn}
Consider the following toy example with three sites ($n=3$),
\begin{align*}
S_1 & = \{ {\tt A123XYZ, C789ABC} \} \\
S_2 & = \{ {\tt A123XYZ, A123XDR, D555SDD} \} \\
S_3 & = \{ {\tt C789ABC, A123XYZ} \}.
\end{align*}
In passing, it should be noted that a formal requirement for something to be a set is that its
members are distinct.  If this formal requirement is not met then each member of the set
could be tagged by a unique number which is not considered in later equality relations.  
This is a technicality which will not be mentioned again and does not affect what follows.

\begin{defn}
A tuple of observations $\y= (y_1, \dots, y_n)$ is an n-tuple consisting of one 
member of each set of observations --- that is, $y_i \in S_i$ for all $i$.
\label{defn:y}
\end{defn}
Continuing the previous example, \\
$\y = ({\tt A123XYZ, A123XYZ, C789ABC})$ \\
is the tuple
formed by taking the first observation from each set.

\begin{defn}
The set of all tuples of observations $\Sc$ in the data is the set of all such $\y$ which can
be formed from the sets $S_1, \dots, S_n$.  This is clearly the cartesian product given by
\begin{equation*}
\Sc = S_1 \times S_2 \dots \times S_n.
\end{equation*}
\label{defn:Sc}
\end{defn}
So, in the example framework given before, then the set $\Sc$ has twelve members and is
given by
\begin{align*}
\Sc = & \{ ({\tt A123XYZ, A123XYZ, C789ABC}), \\
& ({\tt A123XYZ, A123XYZ, A123XYZ}), \\
& \quad \dots ({\tt C789ABC, D555SDD,A123XYZ}) \}.
\end{align*}
Considering, the members of $\Sc$ it is obvious that \\
$( {\tt A123XYZ, A123XYZ, A123XYZ} )$\\
is the type of observation which is most of interest, the same individual observed
across all three sites under investigation.  Also, in some way, the
tuples \\
$({\tt A123XYZ, A123XDR, A123XYZ})$ \\and\\ $({\tt C789ABC, A123XDR, C789ABC})$ \\
are
in some way structurally similar (they match at sites one and three) and both are
structurally different to\\ $({\tt A123XYZ, A123XYZ, C789ABC})$.\\
This structural similarity
will now be formalised by using the concept of a {\em type of match}.

\begin{defn}
Two n-tuples of observations $\y= (y_1, \dots y_n)$ and $\z= (z_1, \dots z_n)$ are the 
same type of match ($\y \sim \z$) if whenever two elements of $\y$ match then the same
two elements of $\z$ match and vice versa.  Formally,
\begin{align*}
\y \sim \z \text{ if and only if } & (y_i = y_j) \Leftrightarrow (z_i = z_j) 
\\ & \text{ for all } 
i, j\in \{1, 2, \dots, n\}. 
\end{align*}
\label{defn:sim}
\end{defn}
Note that, for simplicity the limits
$i, j \in \{1, 2, \dots, n\}$ on indices
will usually be omitted where, as in this case, they are obvious.
It can trivially be shown that the relation defined by $\sim$ meets the requirements
of an equivalence relation in set theory.

\section{The set of every type of match}

Having formalised the concept of when two sets of observations are the same type of
match, the next step is to introduce an entity which can represent the type of match
of a given tuple of observations.  This is simply achieved using partitions of the
first $n$ integers.  A partition of the first $n$ integers is a set of sets 
$\Mb = \{ M_1, M_2, \dots M_m \}$ such that each integer from one to $n$ is in one and
only one of the sets $M_1 \dots M_m$.  (In the literature, these $M_i$ are often
referred to as blocks.)  Any n-tuple of observations is
related to some such $\Mb$ by the relation given in Definition \ref{defn:Cy}.

\begin{defn}
A {\em type of match} is a partition $\Mb$ of the first $n$ integers which is used to
represent the structure of matches within an n-tuple of observations $\y$.  The 
relationship between $\Mb$ and $\y$ is given by Definition \ref{defn:Cy}.
\label{defn:Mb}
\end{defn}
Considering the first three integers, then $\{ \{1,2,3\} \}$, $\{ \{1,2\}, \{3\} \}$
and $\{ \{ 1 \}, \{ 2 \}, \{ 3 \} \}$ are among the possible partitions.

\begin{defn}
The set $\Mn$ is the set of all possible partitions of the first $n$ integers.  This
can be used to represent any possible type of match over $n$ observation sites.
\label{defn:Mn}
\end{defn}
For one site only the partition $\{ \{ 1 \} \}$ is in $\Mc_1$.  For two sites,
two possible partitions are available $\{ \{1,2\} \}$ and $\{ \{1\}, \{2\} \}$.  For
three sites, five partitions are avaialble.  The enumeration of $\# \Mn$ is well
understood and uses the Bell numbers \cite{biggs61}. The sequence of the Bell
numbers begins 1, 2, 5, 15, 52, 203, 877, 4140, 21147.

\begin{defn}
The type of match of an n-tuple of observations $\y= (y_1, \dots, y_n)$ is
given by $C(\y) \in \Mn$ where $C(\y)= \Mb = \{ M_1, \dots, M_m\}$ is the partition
of the first $n$ integers which satisfies $(y_i = y_j) \Leftrightarrow i, j \in M_k$ 
for some $k \in [1, 2, \dots, m]$.  That is, $\Mb$ is the partition chosen such that
any two site indices are in the same block within $\Mb$ if and only if the observations in $\y$
at those sites are equal.
\label{defn:Cy}
\end{defn}
It can clearly be seen that $C(\y)$ is uniquely specified by this definition.
To continue with the earlier example, if \\
$\y = ({\tt A123XYZ, A123XYZ, C789ABC})$ then\\
$C(\y) = \{ \{1,2\} \{3\} \}$ \\and if \\
$\y = ({\tt A123XYZ, A123XYZ, A123XYZ})$ then\\
$C(\y) = \{ \{1,2,3\} \}$.

It must now be shown that $C(\y)$ works as a representation of the type of match in
a consistent way with the relationship $\sim$ given by Definition \ref{defn:sim}.

\begin{thm}
For n-tuples of observations $\y = (y_1, \dots, y_n)$ and $\z = (z_1, \dots, z_n)$ then
\begin{equation*}
C(\y) = C(\z) \text{ if and only if } \y \sim \z.
\end{equation*}
\end{thm}
\begin{pf}
Let $\Mb_y = C(\y)$ and $\Mb_z = C(\z)$.
First it must be shown that $(\y \sim \z) \Rightarrow (\Mb_y = \Mb_z)$.  This follows
trivially.  Since $(y_i = y_j) \Leftrightarrow (z_i = z_j)$ then if $i,j$ are in the same
set in $\Mb_y$ they must be in the same set in $\Mb_z$ and if they are in different
sets in $\Mb_y$ they must be in different sets in $\Mb_z$.  As all integers from one to $n$
appear once each in both partitions then it must be the case that $\Mb_y = \Mb_z$.

Similarly it must be shown that $(\Mb_y = \Mb_z) \Rightarrow (\y \sim \z)$.  A very similar
argument applies.  If $i,j$ are in the same set in $\Mb_y$ (and therefore in $\Mb_z$)
then $y_i = y_j$ and also
$z_i = z_j$ if they are in different sets then $y_i \neq y_j$ and also $z_i \neq z_j$.  
Therefore $(y_i = y_j) \Leftrightarrow (z_i = z_j)$ and hence $\y \sim \z$.
\end{pf}
It is useful at this point to define a shorthand notation for the type 
of match of most interest,
that where the observations are the same at every site.
\begin{defn}
Let $\Ab_n \in \Mn$ represent a {\em true match}, that is the type of match where the same
observation is made over all $n$ sites.  Therefore,
$\Ab_n = \{ \{1,2,\dots, n \} \}$.
\label{defn:taun}
\end{defn}

\section{Introducing false matching into the framework}

So far the false match problem has been ignored and it has been assumed
that for a given n-tuple of observations $\y= (y_1,\dots, y_n)$ then the relation
$y_i = y_j$ can be taken at face value.  However, the original problem was that,
in licence plates, partial observations can lead to two distinct individuals
being confused.  In order to capture this in the described framework, a partial
ordering  
will be introduced on the set $\Mn$ and this will then be related
to the partial observations.  (It is somewhat unfortunate that this paper uses the
phrase ``partial plate survey'' from transportation and the 
term ``partial ordering'' from set theory.  These terms should not be confused.)

The next step is to introduce a partial ordering on the set $\Mn$.  It will be seen
in the next section how this relates to the false matching problem.
\begin{defn}
For two partitions \\
$\Mb = \{M_1,\dots, M_m\} \in \Mn$ \\
and \\
$\Mb' = \{ M'_1, \dots, M'_{m'} \} \in \Mn$\\
a partial ordering $\succsim$ is given by,
\begin{equation*}
\Mb \succsim \Mb' \text{ if and only if } (i,j \in M_k) \Rightarrow (i,j \in M'_l),
\end{equation*}
for some $k$ and $l$.  Put more simply, $\Mb \succsim \Mb'$ 
if whenever $i$ and $j$ are in the same set within $\Mb$ then
they are also in the same set within $\Mb'$.   

The symbol $\succ$ will be used to mean {\em strictly succeeds}.  That is
$\x \succ \y$ means $\x \succsim \y$ and $\x \not \sim \y$. The symbol
$\succ \succ$ will be used to mean {\em immediate successor} that is,
if $\x \succ \succ \z$ then $\x \succ \z$ but there is no $\y$ such
that $\x \succ \y \succ \z$. The symbols $\succ$, $\precsim$ and 
$\prec \prec$ will have their obvious meanings.  
\end{defn}
It can be trivially shown that this relation meets the formal requirements for a
partial ordering.  It should also be noted that this relation is extremely close
to the original equivalence relation but with the implication going in one direction
only.  It can also be shown that under this partial ordering then $\#\Mb$
the number of sets (blocks) in $\Mb \in \Mn$ is a consistent enumeration of $\Mn$.

\begin{figure*}
\begin{center}
\includegraphics[width=14cm]{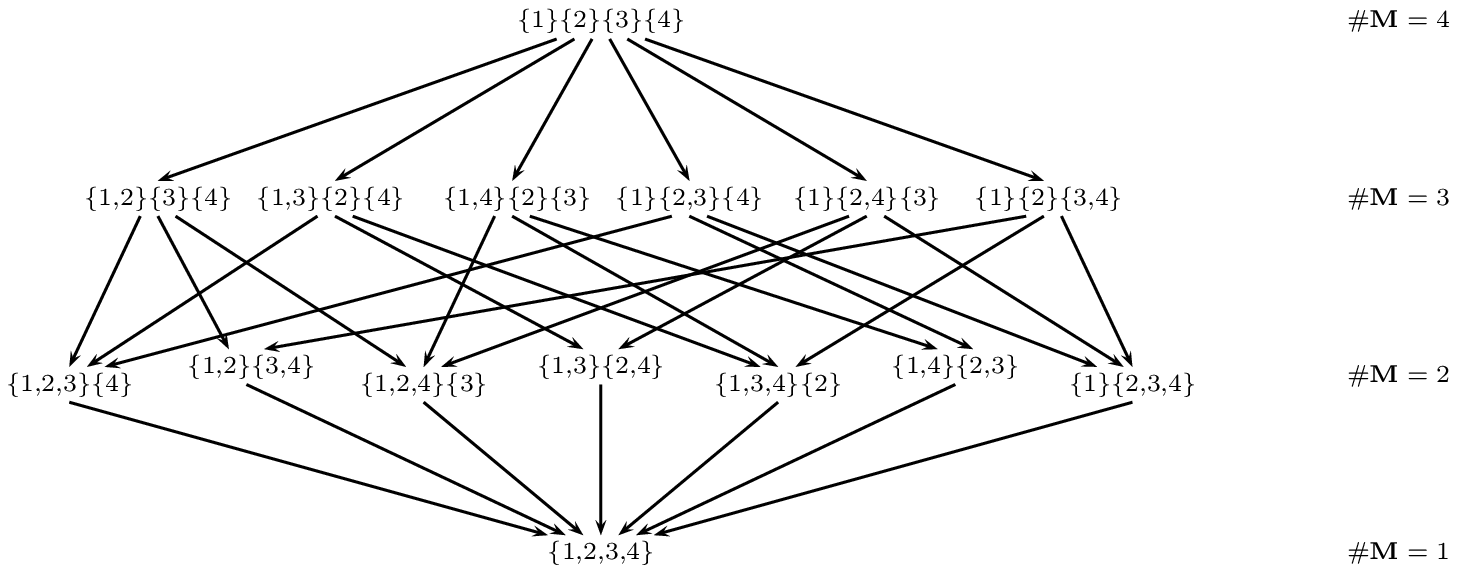}
\end{center}
\caption{Hasse diagram for $\Mc_4$.}
\label{fig:hassem4}
\end{figure*}

A Hasse diagram is a way of visualising a partially ordered set.  A Hasse
diagram is constructed by plotting a partially ordered 
set $S$ graphically in such a way that for all $\x, \y \in S$ if 
$\x \prec \y$ 
then $\x$ is further to the bottom of the diagram than $\y$.
An arrow is drawn in a Hasse
diagram from $x$ to $y$ if $x \succ \succ y$.  Figure
\ref{fig:hassem4} shows the Hasse diagram of $\Mc_4$ with the
partial ordering given by the previous definition.  

\begin{defn}
Given an n-tuple of observations $\y = (y_1, \dots, y_n)$, let 
$\y^* = (y_1^*, \dots, y_n^*)$
represent the {\em partial observation} formed from $\y$.  Since a partial observation can
cause distinct individuals to appear the same but cannot cause the same 
individual to appear distinct at different sites then the following relation
holds,
\begin{equation*}
(y_i = y_j) \Rightarrow (y_i^* = y_j^*).
\end{equation*}
This star notation will also be used to distinguish the set of all possible
partial observations in the data $\Sc^*$ and, in general, to distinguish
functions which apply to partial data rather than the full data.
\label{defn:partial}
\end{defn}
Note that this is the only assumption so far made about the nature of the partial
observation.  In licence plate surveys then the choosing of which part of a plate
to survey needs to be made with reference to the particular format of plate to be
observed.  
Consider the observations from the earlier example.  If 
$\y = ( {\tt A123XYZ, A123XDR, C789ABC} )$ then a standard way to make
partial observations on this type of plate is to collect only the first letter
and the digits.  Therefore $\y^* =  ( {\tt A123, A123, C789} )$.  Note that
$C(\y) \neq C(\y^*)$ since $y_1 \neq y_2$ but $y_1^* = y_2^*$.  The way that
$C(\y)$ can change when only a partial observation is made is given by the next
theorem.

\begin{thm}
If $\y = (y_1, \dots, y_n)$ is an n-tuple of observations then
\begin{equation*}
C(\y^*) \precsim C(\y).
\end{equation*}
\end{thm}
\begin{pf}
Let $\Mb = (M_1, \dots, M_m) = C(\y)$ and $\Mb' = (M'_1, \dots, M'_{m'}) = C(\y^*)$.
The theorem follows trivially from the relation given in Definition
\ref{defn:partial}.
If $i,j \in M_k$ for some $k$ then $y_i = y_j$ and
hence $y^*_i = y^*_j$ which in turn implies, $i,j \in M'_l$ for some $l$.  Therefore,
$(i,j \in M_k) \Rightarrow (i,j \in M'_l)$ which is the condition for the partial
ordering.
\end{pf}

From this theorem, it can be seen that when only partial data is available, the 
type of match of the partial observation may change only in a given way.  Specifically, the
type of match  of the partial data can be the same as that of the full data or any
type of match available by following down the arrows on the Hasse diagram.

Next, some counting functions are defined --- these are used to enumerate the
number of matches in the data which are different types of match.

\begin{defn}
Let $\y$ be an n-tuple of observations and $\Mb \in \Mn$ be a type of match.
The {\em exact matching function} for an observation $\y$ is 
defined by,
\begin{equation*}
x(\y,\Mb) = \begin{cases}
1 & \text{ if and only if } C(\y) = \Mb \\
0 & \text{ otherwise }.
\end{cases}
\end{equation*}
\label{defn:xyM}
\end{defn}

\begin{defn}
Let $\Mb \in \Mn$ be a type of match.
The {\em exact matching function} for $\Sc$ the set of all observations is given by,
\begin{equation*}
X(\Mb) = \sum_{\y \in \Sc} x(\y, \Mb).
\end{equation*}
\label{defn:XM}
\end{defn}
It can be readily seen that $X(\Mb)$ is the number of n-tuples $\y \in \Sc$ which have
a type of match $C(\y) = \Mb$.  It can be further seen that the original problem
of counting the number of individuals seen at all of $n$ sites is the problem of
evaluating $X(\Ab_n)$.

\begin{defn}
Let $\y$ be an n-tuple of observations and $\Mb \in \Mn$ be a type of match.
The {\em relaxed matching function} for an obervation is defined by,
\begin{equation*}
r(\y,\Mb) = \begin{cases}
1 & \text{ if and only if }  C(\y) \precsim \Mb \\
0 & \text{ otherwise }.
\end{cases}
\end{equation*}
Equivalently,
\begin{equation*}
r(\y,\Mb)= \sum_{\Mb' \precsim \Mb} x(\y, \Mb').
\end{equation*}
\label{defn:ryM}
\end{defn}

\begin{defn}
Let $\Mb \in \Mn$ be a type of match.
The {\em relaxed matching function} for $\Sc$ the set of all observations is given by
\begin{equation*}
R(\Mb) = \sum_{\y \in \Sc} r(\y, \Mb).
\end{equation*}
Equivalently,
\begin{equation*}
R(\Mb)= \sum_{\Mb' \precsim \Mb} X(\Mb').
\end{equation*}
\label{defn:RM}
\end{defn}
It should be noted in passing that $R(\Ab_n) = X(\Ab_n)$ since there are no 
$\Mb \prec \Ab_n$.

\section{Solving the false match problem}

In order to solve the false match problem, it is necessary to prove some simple 
lemmas which relate these counting functions.  The main goal here is 
to estimate $X(\Ab_n)$ (the number of n-tuples
representing the same individual at all $n$ sites) in terms of the partial data
$\Sc^*$.  The second goal is to do this in a way which does not involve investigating
every single possible n-tuple.  The reason for this is that a realistic size for
a traffic survey is of the order of one thousand vehicles.  If there are six
sites, then there are $1000^6$ tuples to investigate and this would be far too
slow computationally.

\begin{lem}
Any exact matching function can be expressed in terms of relaxed matching functions and
``lower'' exact matching functions.
\begin{equation*}
X(\Mb) = R (\Mb) - \sum_{\Mb' \prec \Mb} X(\Mb').
\end{equation*}
\label{lem:rel1}
\end{lem}
\begin{pf}
This follows trivially from Definition \ref{defn:RM}.
\end{pf}
This expression can be used recursively so that any
$X(\Mb)$ can be expressed as a function of $R(\Mb')$ for all 
$\Mb' \precsim \Mb$.  The lemma can be thought of as being a version of
the inclusion/exclusion principle for partitions of the integers under this
partial ordering.

\begin{defn}
Let $M= \{m_1, \dots, m_l\} $ be a set of integers, such 
that $m_i \in \{1,2,\dots, n\}$ for all $i$.  Let $\Sc'$ be the set of l-tuples
of observations formed by the cartesian product,
\begin{equation*}
\Sc'= S_{m_1} \times S_{m_2} \times \dots \times S_{m_l}.
\end{equation*}
In other words, $\Sc'$ is the set of l-tuples of observations over some subset of the original 
sites.  Then define,
\begin{equation*}
T(M) = X(\Ab_l),
\end{equation*}
where the exact match $X(\Ab_l)$ is in this case over the l-tuples in $\Sc'$ rather than
the n-tuples in $\Sc$.  In other words, $T(M)$ is the number of individuals seen at all
sites in the set $M$.
\label{defn:TM}
\end{defn}
Note that,
It can be easily seen that the problem of evaluating $T(M)$ is either exactly the same
as the original problem, if $M = \{1, 2, \dots, n\}$ or it is a sub problem over a reduced
number of sites.  If $M$ has a single member $M = \{m\}$ then $T(M)$ is simply the number
of observations in set $S_{m}$ that is, $T(\{m\}) = \#S_{m}$.

\begin{lem}
The relaxed matching function $R(\Mb)$ where $\Mb = \{M_1, \dots, M_m\} \in \Mn$
can be expressed as a product of exact matches over subsets of sites using the 
expression,
\begin{equation*}
R(\Mb)= \prod_{i=1}^m T(M_i).
\end{equation*}
\label{lem:rel2}
\end{lem}
\begin{pf}
Clearly, for an n-tuple of observations $\y = (y_1, \dots, y_n)$ then,
\begin{equation*}
r(\y, \Mb) = \begin{cases} 
1 & \text{ if for all $i,j,k$ then }  \\ 
& (i, j \in M_k) \Rightarrow (y_i = y_j) \\
0 & \text{ otherwise.} 
\end{cases}
\end{equation*}
Therefore,
\begin{align*}
& \sum_{\y \in \Sc} r (\y, \Mb)   = \\
& \# \{\y \in \Sc: (i, j \in M_k) \Rightarrow (y_i = y_j) 
\text{ for all } i,j,k \}.
\end{align*}
The left hand side of this is simply $R(\Mb)$ as required.  Since $\Sc$ is the
cartesian product then it can be seen that those $\y \in \Sc$ which meet the condition
are those which are picked out by $T(M_i)$ and therefore the right hand side is
$\prod_{i=1}^m T(M_i)$ as required.
\end{pf}
Note that if $\Mb=\Ab_n$ then this expression simply says 
$R(\Ab_n) = T(\{1, 2, \dots,n\}) = X(\Ab_n)$.  In all other cases, this allows a
relaxed matching function to be expressed as a product of exact matching functions over a
subset of the original sites.
 
\begin{defn}
The probability $p(i)$ where $i \in \{1,2,\dots,n\}$ is the probability that, 
given that $i$ observed individuals are all different in
the full observation, they will all be the same in the partial observation.  For the
method described to work, this
$p(i)$ must be independent of the sites at which the vehicles are observed.
By convention, $p(1) = 1$.
\label{defn:pi}
\end{defn}
It should be noted that this definition does place some restrictions on the type of
data which can be analysed by this method and which types of partial observations
are suitable.  A discussion of $p(i)$ in the context of licence plate observations
follows this section.  It is likely that other formulations of this problem would
be possible if $p(i)$ varies with the sites considered.
\begin{lem}
An unbiased estimator $\hat{t}$ for $X(\Ab_n)$ is given by,
\begin{equation*}
\hat{t} = X^*(\Ab_n) - \sum_{\Mb \succ \Ab_n} p(\#\Mb)X(\Mb).
\end{equation*}
\label{lem:est}
\end{lem}
\begin{pf}
The quantity $X^*(\Ab_n)$ is equal to $X(\Ab_n)$ plus all those n-tuples of
observations which are false matches.  Each element of the sum represents the 
number of false matches arising from a given type of match.  Writing this
out formally,
\begin{align*}
\hat{t} = & X^*(\Ab_n) - \\ \sum_{\Mb \succ \Ab_n} 
& \E{ \# \{\y \in \Sc : C(\y^*) = \Ab_n, C(\y) = \Mb\} }.
\end{align*}
The set $\{\y \in \Sc : C(\y^*) = \Ab_n, C(\y) = \Mb\}$ is the set of
n-tuples in the data $\Sc$ which are a match of type $\Mb$ in the complete
data but appear to be a match of type $\Ab_n$ in the partial data $\Sc^*$.
Now, the number of distinct individuals
in this n-tuple must be $\#\Mb$.  Therefore,
\begin{equation*}
\Prob{C(\y^*) = \Ab_n | C(\y) = \Mb) }= p(\#\Mb).
\end{equation*}
Bayes theorem gives,
\begin{align*}
&\Prob{C(\y^*) = \Ab_n , C(\y) = \Mb)} \\
& =  p(\#\Mb) \Prob{C(\y) = \Mb} \\
& = \frac{p(\#\Mb) X(\Mb)}{\#\Sc}.
\end{align*}
Hence, the expected number of false matches arising from each type of
match can be given by,
\begin{align*}
& \E{ \#\{\y \in \Sc : C(\y^*) = \Ab_n, C(\y) = \Mb \} } \\ 
& = \#\Sc\Prob{C(\y^*) = \Ab_n , C(\y) = \Mb)},
\end{align*}
and the lemma follows immediately.
\end{pf}
It may not be immediately obvious that Lemmas \ref{lem:est}, \ref{lem:rel1} and 
\ref{lem:rel2} together allow an unbiased estimate of the number of true matches, 
from the partial plate data (assuming that the $p(i)$ are known).  
First, looking at Lemma \ref{lem:est}, the
quantity $ X^*(\Ab_n)$ can be simply enumerated by computer in the partial data.  
Therefore, this lemma allows an unbiased estimate of the number of matches in the
complete data if an unbiased estimate of $X(\Mb)$ can be found for all 
$\Mb \succ \Ab_n$.  Now, Lemma \ref{lem:rel1} allows $X(\Mb)$ to be expressed
as a sum of $R(\Mb')$ for all $\Mb' \precsim \Mb$.  Lemma \ref{lem:rel2} allows
those $R(\Mb')$ to be either equal to the original required quantity $X(\Ab_n)$
or to be expressed in terms of a product involving subproblems on a reduced number
of sites.  Hence, computer algebra can be used to give an equation which is in terms
of $X(\Ab_n)$ (the quantity desired), $X^*(\Ab_n)$ (measureable on the data), $p(i)$
(assumed to be known) and $T(M)$ (which is a subproblem of the original problem with
a reduced number of sites).  The computer can then be used to recursively solve the
subproblem which has already been shown to be trivial for just one site.  An
expanded description of this solution process is given in \cite[Chapter 4]{clegg2004}.

\subsection{Estimating the probability of false matches}

The method described here relies on a good estimate of $p(i)$ and also on the
assumption that this does not vary by the sites chosen.  The specific details
of British licence plates are not of general interest (and it should again be stressed
that the method discussed here is general and not limited just to specific types of
licence plate survey, indeed it could be used for any type of data collection where
the restrictions on $p(i)$ are met).  However, illustrating how $p(i)$ can be
estimated in a practical case might be of interest and illuminate how the method
was applied in real life.  More details on this can be found in \cite[Chapter 5]{clegg2004}.

Two methods of estimating $p(i)$ are practical.  If the distribution of the vehicle types
can be calculated then an analytical approach is possible.  Let there be $N$ vehicle
types which are distinguishable in the partial observations and let $f_j$ be the proportion 
of the vehicle fleet which is of type $j$ (assume 
the membership of each type is relatively large).  Therefore, $p(i)$ is approximately
given by,
\begin{equation*}
p(i) = \sum_{j=1}^N f_j^i.
\end{equation*}
In the case of the old style British licence plates discussed, the distribution of the digits
is almost a flat distribution from 1 to 999.  The distribution of the year letters is more
complex and can be estimated from consideration of the data.  Therefore, $f_j$ can be
calculated for each possible partial observation and hence $p(i)$.

An alternative method is to estimate $p(2)$ by finding two sites which
are so far separated geographically that no vehicle could be seen at both.  Any vehicle
seen at both must be a false match and therefore if there are $x$ observed matches in
the partial data and then $\hat{p(2)} = x/(\#S_1 \#S_2)$.  Similarly $p(3)$ can be
estimated by finding three such geographically distant sites.  Higher order $p(i)$ can
be estimated with reference to the previous method or by assuming a functional form for the
fall off.  An estimate of $p(2) = 7.4 \times 10^{-6}$ was given
in \cite[Chapter 5]{clegg2004} for licence plate data of the type discussed.

\section{Results on simulated data}

Table \ref{4tab:simres} shows simulation results for between two and
six observation sites. These could be thought of as one site observed on several
days, or six sites observed on several different days.  
Num. Veh. refers to the total number of observations at each of the
sites (in these simulations, there are the same number of vehicles in
each data set). The five columns of the form 
1 -- $n$ refer to the number of vehicles which genuinely
went from site one to site $n$ visiting all sites in between.  If
this column is blank it means that there was no site $n$.
For example, if 1 -- 2 = 100, 1 -- 3 = 200 and
1 -- 4 is blank. This means that 100 vehicles travelled between
site one and site two,  200 vehicles travelled between sites one,
two and three and there were only three sites.
Note that these are cumulative so that if 1 -- 2 = 20
and 1 -- 3 = 10 this means that 30 vehicles in total went from site
one to site two and ten of them continued to site three.
Thus the first experiment is two sites, 1000 vehicles at each for
which there were ten vehicles which were genuinely seen at both sites.
In every experiment, the number of different vehicle types
was set at 10,000 with a flat distribution. 
Note that the simplifying assumptions of a flat distribution
and the same number of vehicles at each site are simply there to make
the experiment easier to understand rather than being necessary for the
method to work.
It should be clear that the desired answer from
the correction process is the rightmost figure in these columns.

Each experiment is repeated twenty times with simulated data being generated
anew each time.  The correction process has no random element and will always
give the same result for the same data.  The mean raw
number of matches is given --- this is the total number of
n-tuples which were seen to have the same value for each observation at
every site (averaged over the twenty simulation runs).
Because of the combinatorial nature of the procedure,
this could, in principle, be much larger than the number of vehicles in any
of the data sets (since it counts any n-tuple).
The sample standard deviation ($\sigma$) is given for
the raw matches.  The mean estimated correct number of matches is then given
(again averaged over the twenty simulations).  The sample standard deviation
$\sigma$ is then given for the twenty corrected matches.  It is clear that
the most important test is that the mean corrected number of matches is
as near to correct as possible.  However, it should also be kept in mind that
in reality, a researcher could only run the matching procedure once on any
given set of data so it is also important that $\sigma$ is as low
as possible.  A significant improvement to the method would be to estimate
the variance as well as producing the mean in order that the researcher could have
some idea as to the likely accuracy of the corrected results. 

The first five rows are all results on just two test sites.  This procedure
is not the ideal one to use for estimates on matches between just two sites
and the work of other authors in the field should be used in such
a circumstance.  However, these results are included here for completeness.
In the first
experiment, the average number of raw matches over the twenty runs is 111.4.  The
average number of corrected matches is 100 less than this (11.4).  This is
close to the correct answer of 10.  However, it should be noticed that
the $\sigma$ is high in comparison to the actual answer.  In this case, the
$\sigma$ is 8.5 which is of the same order of magnitude as the answer.
This is to be expected since we are looking for only 10 true matches in over
110 observed matches.  If we increase the number of vehicles to 2000 then,
as would be expected, the number of false matches goes up (to approximately
400) and the $\sigma$ also rises (to almost 20).

The next five rows of results are all over three sites.  In the first of
these, 10 vehicles travel between all three and all other matches are
coincidence.  1000 vehicles are observed at all sites.  The mean
corrected match across all sites 9.3 is close to the actual answer of 10 and
the $\sigma$ is lower than in the two site case.  However, when the same
experiment is run with 500 vehicles travelling from sites one to two in addition
to 10 vehicles travelling from sites two to three, the $\sigma$ increases markedly
(it almost doubles).  In all cases with three sites, the mean is a good estimate
and the $\sigma$ is generally low enough that a good estimate can be expected.

The next four rows of results are for experiments made over four sites.  The
first experiment has 100 vehicles which visit all four.  The mean corrected match
is 104 (very close) and the $\sigma$ is only 22.  It is hard to explain why
this $\sigma$ actually falls in the next experiment when more vehicles are 
genuinely seen in common between the other sites.  
In all cases the mean of the predictions is approximately correct (the
worst performance being in the case of the fourth experiment when the mean was 106.1
not 100).  

The next six rows of results are experiments made over five sites.  Again, the
mean corrected results are approximately correct.  However, in the worst case,
the mean is 11 too high and the $\sigma$ in the results is 46.7 which is comparable
to the level of the effect being observed.  In this case approximately 120 false matches
are being removed each time.  However, previous experiments have been able to
correct for a greater proportion of false matches with less $\sigma$ in the
result.

The final four rows of results are experiments over six sites.  This was the
largest number of sites for which it was practical to do runs of twenty or more
simulations with the computer power available.  Again, the mean corrected estimate
of matches was nearly correct in all cases.  The worst performance was an
estimate of 92.2 (correct result 100).  The $\sigma$ was, however,
relatively high.  This was a surprise in some cases --- particularly the first row
of results where the mean number of false matches was only 21.2.  In many senses,
the worst results was the final one where a $\sigma$ of 55.0 was given on an
corrected prediction of only 101.3.

The time taken to do one run over six sites with one thousand pieces of data on each site
was thirty seconds on a Celeron 366 computer running Debian Linux.   Six sites 
with one thousand vehicles at each is a reasonable
size for a typical traffic survey.  It is practical
(if time consuming) to do experiments on seven sites, even using such comparatively
obsolete equipment.  However, eight sites or more is probably too computationally
expensive for the moment and this is a limitation of the method outlined.  The exact
rate at which the computational requirements increase with the number of sites
is hard to determine.  It
will relate to the Bell numbers, to the number of observations at each site and to
the number of pairs of observations at each site pair.

The results given here are certainly consistent with the idea that the method
gives an unbiased estimator for the true number of matches.  In some experiments,
there were problems with the standard deviation being higher than would be
desirable in real cases.  It is important to bear in mind that these were relatively
extreme tests of the method since $p(2)$ and $p(3)$ were relatively low and the
number of samples given were quite high.  Often the method was attempting to predict
only ten true matches in a number of observed matches which might be several hundred.

\begin{table*}
\begin{center}
\begin{tabular} { |l| l l l l l|l l l l | } \hline

No.  & 1 -- 2 & 1 -- 3 &
1 -- 4 & 1 -- 5 & 1 -- 6 & Av. Raw &
$\sigma$ Raw & Av. Cor. & $\sigma$ Cor. \\
Veh. & & & & & & Matches & Matches & Matches & Matches \\ \hline

1000 & 10 &  &   &   &   &  111.4 &   8.5 &  11.4 &   8.5 \\
2000 & 10 &  &   &   &  & 411.8 &  19.5 &  11.8 &  19.5 \\
1000 & 100 &   &   &   &   &  199.2 &  12.0 &  99.2 &  12.0 \\
1000 & 200 &   &  &  &   &  302.3 &   7.7 & 202.3 &   7.7 \\
1000 & 500 &  &   &   &   &   596.6 &  12.3 & 496.7 &  12.3 \\  \hline
1000 & 0 & 10 &   &   &   &   21.9 &   4.6 &   9.3 &   3.3 \\
1000 & 500 & 10 &   &   &   &   73.8 &   7.5 &  10.2 &   6.2 \\
1000 & 100 & 100 &   &   &    &  152.1 &   8.5 & 101.9 &   7.5 \\
1000 & 500 & 250 &   &   &   &  388.3 &  22.7 & 253.2 &  20.1 \\
1000 & 0 & 500 &   &  &   &  667.2 &  24.9 & 506.0 &  22.3 \\ \hline
1000 & 0 & 0 & 100  &  &  & 154.6 &  26.6 & 104.0 &  22.6  \\
1000 & 100 & 100 & 100 &    &   &  164.4 &  11.4 &  97.7 &   9.3\\
500 & 100 & 100 & 100  & &  &  140.7 &  19.3 & 105.8 &  17.4 \\
1000 & 500 & 250  & 100  &  &   & 207.8 &  29.7 & 106.1 &  23.7 \\ \hline
500 & 10  & 10 & 10 & 10 &  &  14.2 &   2.2 &  10.5 &   1.8 \\
1000 & 10 & 10 & 10 & 10 &  &  17.4 &   4.1 &   9.4 &   2.8 \\
500 & 50 & 50 & 50 & 50 &  &  71.3 &  14.3 &  47.8 &  12.3 \\
500 & 100 & 100 & 100 & 100 &  &  151.9 &  26.9 &  92.0 &  22.3 \\
1000 & 0 & 0 & 0 & 100 &  &  177.6 &  29.9 & 103.4 &  22.6 \\
1000 & 100 & 100 & 100 & 100 &  & 222.2 &  61.5 & 111.0 &  46.7 \\ \hline
1000 & 0 & 0 & 0 & 0 & 10 &  21.2 &  13.4 &  12.3 &   9.9  \\
500  & 0 & 0 & 0 & 0 & 100 & 152.6 &  45.5 &  92.2 &  37.3 \\
1000 & 0 & 0 & 0 & 0 & 100 & 214.6 &  58.0 & 103.5 &  40.2  \\
1000 & 100 & 100 & 100 & 100 & 100 & 289.8 &  88.4 & 101.3 &  55.0 \\ \hline
\end{tabular}

\caption{Simulation results --- all performed over twenty runs with 10,000
distinct vehicle types.}
\label{4tab:simres}
\end{center}
\end{table*}

To test the method more fully, four very extreme tests were given.  Each of
these tests involved six sites at each of which one thousand vehicles were
observed.  Interacting flows were chosen to cause a large
number of false matches in a diversity of ways.  Because these experiments were
chosen to cause a large number of false matches then one thousand runs of
each experiment were performed.  The averaged results are shown in Table
\ref{table:simres2}.

\begin{table*}
\begin{center}
\begin{tabular} { |l l|l l l l | } \hline

 Experiment & Expected & Av.Raw &
$\sigma$ Raw & Av.Cor. & $\sigma$ Cor. \\
Number & Answer & Matches & Matches & Matches & Matches \\ \hline

1 & 0 & 739 & 305 & 11.9 & 196 \\
2 & 0  & 110 & 45.5 & -0.950 & 27.1 \\
3 & 250  & 836 & 287 & 249 &  205 \\
4 & 500 & 1920 & 531 & 496 & 356 \\ \hline

\end{tabular}
\caption{Simulation results --- all performed over one thousand runs with 10,000
distinct vehicle types.}
\label{table:simres2}
\end{center}
\end{table*}

In experiment one, five hundred vehicles travelled from one to five and five hundred
from two to six.  The remaining five hundred vehicles at sites one and six were
appeared nowhere else.  No vehicles made the complete journey.  As can be seen, 
on average over seven hundred false matches were seen and the standard deviation
between runs was extremely large.  However, the mean was within twelve of the
correct answer (zero) although the standard deviation was large.  In such
extreme circumstances, a single experiment would be next to useless but it
is good evidence that the method was unbiased.

In experiment two, five hundred vehicles travelled from one to three.  Five hundred
vehicles travelled from four to six.  Five hundred vehicles visited only odd numbered
sites and five hundred vehicles visited only even numbered sites.  In this
experiment the corrected mean result was almost exact (within one) and the
standard deviation was much lower than the other three experiments.

In experiment three, two hundred and fifty vehicles travelled to all sites.  Five
hundred vehicles went from site one to three and five hundred from four to six.  
The remaining two hundred and fifty vehicles at each site visited only that single 
site.  As can be seen, the corrected result is almost exactly correct although, 
again, the standard deviation is so high that a single reading would be
worthless.

In experiment four, five hundred vehicles visited every site.  Two hundred and fifty
vehicles went from sites one to three.  Two hundred and fifty vehicles went from
sites four to six.  Two hundred and fifty vehicles visited only sites one and two,
two hundred and fifty vehicles visited only sites three and four and two hundred
and two hundred and fifty vehicles visited only sites five and six.  Again, the
mean of all results is very close (within four vehicles) but the standard deviation
is the highest yet seen.  This is not surprising.  The mean number of raw tuples
of matches averaged nearly 2000, twice the number of vehicles at each site.

These four tests provide a convincing demonstration that the method is, indeed, unbiased
as was shown by theory.

\section{Conclusions}

This paper presented a framework for analysis of surveys where matches are
required over more than two data collection points.  The framework given 
formalises the concept of a type of match using the concept of the equivalence
class.  Further a method is given for evaluating $\Mn$ the set of
all possible types of match
over multiple data sets. 

The framework given is then applied to the problem of false matches --- which
is put into the language of set theory using the concept of a partial ordering.
It is shown how this partial ordering can be used to visualise, by means of
a Hasse diagram, the ways in which false matches can occur in data observed at
multiple sites.  The framework was then used to design and implement an algorithm
which was used to estimate the number of true matches in simulated data.  The
algorithm has also been tested on real data from partial plate surveys. 

This algorithm was implemented and tested on simulated data.  The results show
that the estimator seems to be unbiased and in the majority of
cases tested the standard deviation on the results is low.  The method is suitable
for analysis of matches on
data between three and seven test sites but becomes too computationally intensive
after this point.  A significant improvement to the method would be the estimation
of a variance as well as a corrected number of matches.  A potential weakness of
the method is that it relies on good estimates for $p(i)$.

\bibliography{rgc_ejor04}
\end{document}